

\baselineskip=14pt
\parskip=10pt
\def\halmos{\hbox{\vrule height0.15cm width0.01cm\vbox{\hrule height
  0.01cm width0.2cm \vskip0.15cm \hrule height 0.01cm width0.2cm}\vrule
  height0.15cm width 0.01cm}}
\font\eightrm=cmr8 
\font\eighttt=cmtt8
\magnification=\magstephalf

\def\1{{\overline{1}}}
\def\2{{\overline{2}}}
\parindent=0pt
\overfullrule=0in

\def\frac#1#2{{#1 \over #2}}
\bf
\centerline
{
The Sagan-Savage Lucas-Catalan Polynomials Have Positive Coefficients
}
\rm
\bigskip
\centerline{ {\it
Shalosh B.
EKHAD}\footnote{$^1$}
{\eightrm  \raggedright
Department of Mathematics, Rutgers University (New Brunswick),
Hill Center-Busch Campus, 110 Frelinghuysen Rd., Piscataway,
NJ 08854-8019, USA.
{\eighttt c/o zeilberg  at math dot rutgers dot edu} ,
\hfill \break
{\eighttt http://www.math.rutgers.edu/\~{}zeilberg/ekhad.html} .
Jan. 17, 2011.
Exclusively published in the Personal Journal of Shalosh B. Ekhad and Doron Zeilberger .
}
}

In the last few minutes of Bruce Sagan's wonderful talk[1]
about his joint work with Carla Savage[2] (whose notation I shall use),
he mentioned that they can prove that 
$\frac{1}{\{ n+1 \} } \left \{ { {2n} \atop {n} } \right \}$
are polynomials with integer coefficients, and {\it conjectured} that they are in fact
polynomials with {\it positive} coefficients.

But this follows immediately from the identity
$
\frac{1}{\{ n+1 \} } \left \{ { {2n} \atop {n} } \right \}=\left \{ {{2n-1} \atop {n-1} } \right \}+t\left \{ { {2n-1} \atop {n-2} } \right \}$
that after routine cancellations is equivalent to
$\{ 2n \}= \{n+1 \} \{ n \} + t \{n-1 \} \{ n \}$, that is the case $m=n$ of Lemma 2.1 of [2]. \halmos

{\bf References}

[1] Bruce Sagan, {\it Combinatorial Interpretations of Binomial Coefficient Analogues Related to Lucas Sequences},
talk at the Rutgers University Experimental Mathematics seminar on Dec. 9, 2010,
videotaped by Edinah Gnang.  http://www.youtube.com/watch?v=Fdn890jg2U0 \quad .

[2]  Bruce Sagan and Carla Savage, {\it Combinatorial Interpretations of Binomial Coefficient Analogues Related to Lucas Sequences},
Integers {\bf 10} (2010), 697-703, A52.  http://arxiv.org/abs/0911.3159 \quad .

\end